\documentclass[12pt]{article}

\usepackage{amsmath, epsfig, cite, setspace}
\usepackage{amssymb}
\usepackage{amsfonts}
\usepackage{latexsym}
\usepackage{amsthm}

\makeatletter
\renewcommand{\@seccntformat}[1]{{\csname the#1\endcsname}.\hspace{.5em}}
\makeatother

\newtheorem{thm}{Theorem}[section]

\newtheorem{lem}[thm]{Lemma}

\renewcommand{\qed}{\hfill$\Box$\medskip}
\renewcommand{\thefootnote}{*}

\numberwithin{equation}{section}

\begin{document}

\begin{center}
{\large\bf On a supercongruence conjecture of Z.-W. Sun}
\end{center}

\vskip 2mm \centerline{Guo-Shuai Mao}
\begin{center}
{\footnotesize $^1$Department of Mathematics, Nanjing
University of Information Science and Technology, Nanjing 210044,  People's Republic of China\\
{\tt maogsmath@163.com  } }
\end{center}


\vskip 0.7cm \noindent{\bf Abstract.}
In this paper, we partly prove a supercongruence conjectured by Z.-W. Sun in 2013. Let $p$ be an odd prime and let $a\in\mathbb{Z}^{+}$. Then if $p\equiv1\pmod3$, we have
\begin{align*}
\sum_{k=0}^{\lfloor\frac{5}6p^a\rfloor}\frac{\binom{2k}k}{16^k}\equiv\left(\frac{3}{p^a}\right)\pmod{p^2},
\end{align*}
where $\left(\frac{\cdot}{\cdot}\right)$ is the Jacobi symbol.



\renewcommand{\thefootnote}{**}

\section{Introduction}
In the past years, congruences for sums of binomial coefficients have attracted the attention of many researchers (see, for instance, \cite{az-amm-2017,G,gl-arxiv-2019,mc-crmath-2019,mw-ijnt-2019,st-aam-2010,st-ijnt-2011,ws-ijnt-2019}). In 2011, Sun \cite{st-ijnt-2011} proved that for any odd prime $p$ and $a\in\mathbb{Z}^{+}$,
\begin{align*}
\sum_{k=0}^{p^a-1}\binom{2k}k\equiv\left(\frac{p^a}3\right)\pmod{p^2}.
\end{align*}
Liu and Petrov \cite{lp-aam-2020} showed some congruences on sums of $q$-binomial coefficients.

In 2006, Adamchuk \cite{adamchuk-oeis-2006} conjectured that for any prime $p\equiv1\pmod3$,
\begin{align*}
\sum_{k=1}^{\frac{2}3(p-1)}\binom{2k}k\equiv0\pmod{p^2}.
\end{align*}
Recently, Mao \cite{m-arxiv-2020} confirmed this conjecture.\\
Pan and Sun \cite{ps-scm-2014} proved that for any prime $p\equiv1\pmod 4$ or $1<a\in\mathbb{Z}^{+}$,
$$\sum_{k=0}^{\lfloor\frac34p^a\rfloor}\frac{\binom{2k}{k}}{(-4)^k}\equiv\left(\frac{2}{p^a}\right)\pmod{p^2}.$$
In 2017, Mao and Sun \cite{ms-2019-rama} showed that for any prime $p\equiv1\pmod 4$ or $1<a\in\mathbb{Z}^{+}$,
$$\sum_{k=0}^{\lfloor\frac34p^a\rfloor}\frac{\binom{2k}{k}^2}{(16)^k}\equiv\left(\frac{-1}{p^a}\right)\pmod{p^3}.$$
Sun \cite{sun-tj-2013} proved that for any odd prime $p$ and $a\in\mathbb{Z}^{+}$, we have
\begin{align}\label{st}
\sum_{k=0}^{(p^a-1)/2}\frac{\binom{2k}k}{16^k}\equiv\left(\frac{3}{p^a}\right)\pmod{p^2}.
\end{align}
In this paper, we partly prove Sun's conjecture \cite[Conjecture 1.2(i)]{sun-tj-2013}.
\begin{thm}\label{Thsun} Let $p$ be an odd prime and let $a\in\mathbb{Z}^{+}$. If $p\equiv1\pmod3$, then
$$\sum_{k=0}^{\lfloor\frac56p^a\rfloor}\frac{\binom{2k}k}{16^k}\equiv\left(\frac{3}{p^a}\right)\pmod{p^2}.$$
\end{thm}
We shall prove Theorem \ref{Thsun} in Section 2.
\section{Proof of Theorem \ref{Thsun}}
\begin{lem}\label{L} {\rm (\cite{L})}. For any prime $p>3$, we have the following congruences modulo $p$
\begin{align*}H_{\lfloor p/2\rfloor}\equiv-2q_p(2),\ H_{\lfloor p/3\rfloor}\equiv-\frac32q_p(3),\ H_{\lfloor p/6\rfloor}\equiv-2q_p(2)-\frac 32q_p(3).
\end{align*}
\end{lem}
\noindent{\it Proof of Theorem \ref{Thsun}}. In view of (\ref{st}), we just need to verify that
\begin{align}\label{p2}
\sum_{k={(p^a+1)/2}}^{\lfloor\frac56p^a\rfloor}\frac{\binom{2k}k}{16^k}\equiv0\pmod{p^2}.
\end{align}
Let $k$ and $l$ be positive integers with $k+l=p^a$ and $0<l<p^a/2$. In view of \cite{ps-scm-2014}, we have
\begin{align}\label{l2l}
\frac{l}2\binom{2l}l=\frac{(2l-1)!}{(l-1)!^2}\not\equiv0\pmod{p^a}
\end{align}
and
\begin{align}\label{2k2l}
\binom{2k}k\equiv-p^a\frac{(l-1)!^2}{(2l-1)!}=-\frac{2p^a}{l\binom{2l}l}\pmod{p^2}.
\end{align}
So we have
\begin{align*}
\sum_{k={(p^a+1)/2}}^{\lfloor\frac56p^a\rfloor}\frac{\binom{2k}k}{16^k}\equiv\sum_{k={(p^a+1)/2}}^{\lfloor\frac56p^a\rfloor}\frac{-2p^a}{(p^a-k)\binom{2p^a-2k}{p^a-k}16^k}=\frac{-2p^a}{16^{p^a}}\sum_{k=\lfloor\frac{p^a} 6\rfloor+1}^{(p^a-1)/2}\frac{16^k}{k\binom{2k}{k}}\pmod{p^2}.
\end{align*}
It is easy to see that for $k=1,2,\ldots,(p^a-1)/2$,
\begin{align}\label{pa-12}
\frac{\binom{(p^a-1)/2}k}{\binom{2k}k/(-4)^k}=\frac{\binom{(p^a-1)/2}k}{\binom{1/2}k}=\prod_{j=0}^{k-1}\frac{(p^a-1)/2-j}{-1/2-j}=\prod_{j=0}^{k-1}\left(1-\frac{p^a}{2j+1}\right)\equiv1\pmod p.
\end{align}
This, with Fermat little theorem yields that
$$\sum_{k={(p^a+1)/2}}^{\lfloor\frac56p^a\rfloor}\frac{\binom{2k}k}{16^k}\equiv-\frac{p^a}{8}\sum_{k=\lfloor\frac{p^a}6\rfloor+1}^{(p^a-1)/2}\frac{(-4)^k}{k\binom{(p^a-1)/2}{k}}
\equiv-p^a\sum_{k=\lfloor\frac{p^a}6\rfloor}^{(p^a-3)/2}\frac{(-4)^k}{\binom{(p^a-3)/2}{k}}\pmod{p^2}.$$
Thus, by (\ref{p2}) we only need to show that
\begin{align}\label{pa-3}
p^{a-1}\sum_{k=\lfloor\frac{p^a}6\rfloor}^{(p^a-3)/2}\frac{(-4)^k}{\binom{(p^a-3)/2}{k}}\equiv0\pmod{p}.
\end{align}
Now we set $n=(p^a-1)/2, m=\lfloor\frac{p^a}6\rfloor, \lambda=-4$, then we only need to prove that
\begin{align}\label{n-1m}
p^{a-1}\sum_{k=m}^{n-1}\frac{\lambda^k}{\binom{n-1}k}\equiv0\pmod p.
\end{align}
In view of \cite{SWZ}, we have
$$
\sum_{k=m}^{n-1}\frac{\lambda^k}{\binom{n-1}k}=n\sum_{k=0}^{n-1-m}\frac{\lambda^{m+k}}{(\lambda+1)^{k+1}}\sum_{i=0}^{n-1-m-k}\frac{(-1)^i\binom{n-1-m-k}i}{m+i+1}+\frac{n\lambda^n}{(\lambda+1)^{n+1}}\sum_{k=m}^{n-1}\frac{(\lambda+1)^{k+1}}{k+1}.
$$
It is easy to check that for each $0\leq k\leq n-1-m$
\begin{align*}
\sum_{i=0}^{n-1-m-k}\binom{n-1-m-k}i\frac{(-1)^i}{m+i+1}&=\int_0^1\sum_{i=0}^{n-1-m-k}\binom{n-1-m-k}i(-x)^ix^mdx\\
&=\int_0^1x^m(1-x)^{n-1-m-k}dx=B(m+1,n-m-k),
\end{align*}
where $B(P,Q)$ stands for the beta function. It is well known that the beta function relate to gamma function:
$$B(P,Q)=\frac{\Gamma(P)\Gamma(Q)}{\Gamma(P+Q)}.$$
So
$$
B(m+1,n-m-k)=\frac{\Gamma(m+1)\Gamma(n-m-k)}{\Gamma(n-k+1)}=\frac{m!(n-m-k-1)!}{(n-k)!}=\frac1{(m+1)\binom{n-k}{m+1}}.
$$
Therefore
\begin{align*}
\sum_{k=m}^{n-1}\frac{\lambda^k}{\binom{n-1}k}&=\frac{n}{m+1}\sum_{k=0}^{n-1-m}\frac{\lambda^{m+k}}{(\lambda+1)^{k+1}\binom{n-k}{m+1}}+\frac{n\lambda^n}{(\lambda+1)^{n+1}}\sum_{k=m}^{n-1}\frac{(\lambda+1)^{k+1}}{k+1}\\
&=\frac{n}{m+1}\sum_{k=m+1}^{n}\frac{\lambda^{m+n-k}}{(\lambda+1)^{n-k+1}\binom{k}{m+1}}+\frac{n\lambda^n}{(\lambda+1)^{n+1}}\sum_{k=m+1}^{n}\frac{(\lambda+1)^{k}}{k}\\
&=\frac{n\lambda^n}{(\lambda+1)^{n+1}}\left(\frac{\lambda^m}{m+1}\sum_{k=m+1}^{n}\frac{(\lambda+1)^{k}}{\lambda^{k}\binom{k}{m+1}}+\sum_{k=m+1}^{n}\frac{(\lambda+1)^{k}}{k}\right).
\end{align*}
By (\ref{n-1m}), we just need to show that
\begin{equation}\label{nm+1}
p^{a-1}\frac{\lambda^m}{m+1}\sum_{k=m+1}^{n}\frac{(\lambda+1)^{k}}{\lambda^{k}\binom{k}{m+1}}\equiv-p^{a-1}\sum_{k=m+1}^{n}\frac{(\lambda+1)^{k}}{k}\pmod p.
\end{equation}
It is obvious that
\begin{align*}
\sum_{k=m+1}^{n}\frac{(\lambda+1)^{k}}{\lambda^{k}\binom{k}{m+1}}=\sum_{k=m+1}^n\frac{1}{\binom{k}{m+1}}\left(\frac34\right)^k=\sum_{k=m+1}^n\frac1{\binom{k}{m+1}}\sum_{j=0}^k\frac{\binom{k}j}{(-4)^j}=\mathfrak{B}+\mathfrak{C},
\end{align*}
where
$$
\mathfrak{B}=\sum_{j=m+1}^n\frac1{(-4)^j}\sum_{k=j}^n\frac{\binom{k}j}{\binom{k}{m+1}},\ \ \ \ \ \ \ \mathfrak{C}=\sum_{j=0}^{m}\frac1{(-4)^j}\sum_{k=m+1}^n\frac{\binom{k}j}{\binom{k}{m+1}}.
$$
By the following transformation
$$
\frac{\binom{k}j}{\binom{k}{m+1}}=\frac{k!(m+1)!(k-m-1)!}{j!(k-j)!k!}=\frac{(m+1)!(k-m-1)!(j-m-1)!}{j!(k-j)!(j-m-1)!}=\frac{\binom{k-m-1}{j-m-1}}{\binom{j}{m+1}}.
$$
We have
\begin{align*}
\mathfrak{B}=\sum_{j=m+1}^n\frac1{(-4)^j}\sum_{k=j}^n\frac{\binom{k-m-1}{j-m-1}}{\binom{j}{m+1}}=\sum_{j=m+1}^n\frac{1}{(-4)^j\binom{j}{m+1}}\sum_{k=0}^{n-j}\binom{k+j-m-1}{j-m-1}.
\end{align*}
By \cite[(1.48)]{g-online}, we have
$$
\mathfrak{B}=\sum_{j=m+1}^n\frac{1}{(-4)^j\binom{j}{m+1}}\binom{n-m}{j-m}.
$$
It is easy to show that
$$
\frac{\binom{n-m}{j-m}}{\binom{j}{m+1}}=\frac{(n-m)!(m+1)!(j-m-1)!}{j!(n-j)!(j-m)!}=\frac{n+1}{j-m}\frac{\binom{n}j}{\binom{n+1}{m+1}}.
$$
Thus,
$$\mathfrak{B}=\frac{n+1}{\binom{n+1}{m+1}}\sum_{j=m+1}^n\frac{\binom{n}{j}}{(j-m)(-4)^j}.$$
Now we calculate $\mathfrak{C}$. First we have the following transformation
$$
\frac{\binom{k}j}{\binom{k}{m+1}}=\frac{k!(m+1)!(k-m-1)!}{j!(k-j)!k!}=\frac{(m+1)!(k-m-1)!(m-j+1)!}{j!(k-j)!(m-j+1)!}=\frac{\binom{m+1}{j}}{\binom{k-j}{m-j+1}}.
$$
Thus,
$$
\mathfrak{C}=\sum_{j=0}^m\binom{m+1}j\frac1{(-4)^j}\sum_{k=m+1}^n\frac1{\binom{k-j}{m-j+1}}=\sum_{j=0}^m\binom{m+1}j\frac1{(-4)^j}\sum_{k=0}^{n-m-1}\frac1{\binom{k+m+1-j}{m-j+1}}.
$$
By using package \texttt{Sigma}, we find the following identity,
$$
\sum_{k=0}^N\frac1{\binom{k+i}i}=\frac{i}{i-1}-\frac{N+1}{(i-1)\binom{N+i}N}.
$$
Substituting $N=n-m-1, i=m+1-j$ into the above identity, we have
$$
\mathfrak{C}=\sum_{j=0}^{m-1}\binom{m+1}j\frac1{(-4)^j}\left(\frac{m+1-j}{m-j}-\frac{n-m}{(m-j)\binom{n-j}{n-m-1}}\right)+(m+1)\left(-\frac14\right)^m\sum_{k=1}^{n-m}\frac1k.
$$
It is easy to check that
$$
\frac{(n-m)\binom{m+1}j}{\binom{n-j}{n-m-1}}=\frac{(m+1)!((n-m)!(m+1-j)!}{j!(n-j)!(m+1-j)!}=\frac{(m+1)!((n-m)!}{j!(n-j)!}=\frac{(n+1)\binom{n}j}{\binom{n+1}{m+1}}.
$$
Therefore
$$
\mathfrak{C}=(m+1)\sum_{j=0}^{m-1}\frac{\binom{m}j }{(m-j)(-4)^j}-\frac{n+1}{\binom{n+1}{m+1}}\sum_{j=0}^{m-1}\frac{\binom{n}j }{(m-j)(-4)^j}+(m+1)\left(-\frac14\right)^m\sum_{k=1}^{n-m}\frac1k.
$$
Hence
\begin{equation*}
\mathfrak{B}+\mathfrak{C}=(m+1)\sum_{j=0}^{m-1}\frac{\binom{m}j}{(m-j)(-4)^j}+\frac{n+1}{\binom{n+1}{m+1}}\sum_{\substack{j=0\\j\neq m}}^{n}\frac{\binom{n}j}{(j-m)(-4)^j}+(m+1)\left(-\frac14\right)^m\sum_{k=1}^{n-m}\frac1k.
\end{equation*}
That is
\begin{equation}\label{b+c}
\frac{\lambda^m}{m+1}(\mathfrak{B}+\mathfrak{C})=\lambda^m\sum_{j=0}^{m-1}\frac{\binom{m}j}{(m-j)(-4)^j}+\frac{\lambda^m}{\binom{n}{m}}\sum_{\substack{j=0\\j\neq m}}^{n}\frac{\binom{n}j}{(j-m)(-4)^j}+H_{n-m}.
\end{equation}
It is obvious that
\begin{align*}
&\sum_{k=1}^n\frac{(-3)^k}k=\int_0^1\sum_{k=1}^n(-3)^kx^{k-1}dx=-3\int_0^1\sum_{k=0}^{n-1}(-3x)^kdx=-3\int_0^1\frac{1-(-3x)^n}{1+3x}dx\\
&=3\int_0^1\sum_{k=1}^n\binom{n}k(-1)^k(1+3x)^{k-1}dx=\int_1^4\sum_{k=1}^n(-1)^ky^{k-1}dy=\sum_{k=1}^n\binom{n}k(-1)^k\frac{4^k-1}k
\end{align*}
and
\begin{align*}
\sum_{k=1}^n\binom{n}k\frac{(-1)^k}k&=\int_0^1\sum_{k=1}^n\binom{n}k(-1)^kx^{k-1}dx=\int_0^1\frac{(1-x)^n-1}xdx=\int_0^1\frac{y^n-1}{1-y}dy\\
&=-\int_0^1\sum_{k=0}^{n-1}y^kdy=-\sum_{k=0}^{n-1}\frac1{k+1}=-\sum_{k=1}^{n}\frac1{k}.
\end{align*}
This, with \cite[(1.20)]{st-aam-2010} yields that
$$\sum_{k=1}^n\frac{(\lambda+1)^k}{k}=\sum_{k=1}^n\frac{(-3)^k}k=\sum_{k=1}^n\binom{n}k\frac{(-4)^k}k+H_n.$$
On the other hand, by \cite[(1.48)]{g-online} we have
\begin{align*}
\sum_{k=1}^m\frac{(\lambda+1)^k-1}{k}&=\sum_{k=1}^m\frac{(-3)^k-1}{k}=\sum_{k=1}^m\frac1k\sum_{j=1}^k\binom{k}j(-4)^j=\sum_{j=1}^m\frac{(-4)^j}{j}\sum_{k=j}^m\binom{k-1}{j-1}\\
&=\sum_{j=1}^m\frac{(-4)^j}{j}\binom{m}j=(-4)^m\sum_{j=0}^{m-1}\frac{1}{(m-j)(-4)^j}\binom{m}j.
\end{align*}
Hence
$$
\sum_{k=1}^m\frac{(\lambda+1)^k}{k}=(-4)^m\sum_{j=0}^{m-1}\frac{\binom{m}j}{(m-j)(-4)^j}+H_m.
$$
So
\begin{align}\label{m+1nk}
\sum_{k=m+1}^n\frac{(\lambda+1)^k}{k}=\sum_{k=1}^n\binom{n}k\frac{(-4)^k}k+H_n-\lambda^m\sum_{j=0}^{m-1}\frac{\binom{m}j}{(m-j)(-4)^j}-H_m.
\end{align}
In view of \cite[(1.20)]{st-aam-2010}, and by (\ref{l2l}), (\ref{2k2l}) and (\ref{pa-12}) we have
\begin{align}\label{pa-1n}
p^{a-1}\sum_{k=1}^n\binom{n}k\frac{(-4)^k}k\equiv p^{a-1}\sum_{k=1}^n\frac{\binom{2k}k}k\equiv p^{a-1}\sum_{k=1}^{p^a-1}\frac{\binom{2k}k}k\equiv0\pmod p.
\end{align}
It is obvious that
$$
p^{a-1}H_n=p^{a-1}\sum_{k=1}^n\frac1k\equiv p^{a-1}\sum_{j=1}^{(p-1)/2}\frac1{jp^{a-1}}=H_{(p-1)/2}\pmod p
$$
and $p^{a-1}H_m\equiv H_{\lfloor p/6\rfloor}\pmod p$, $p^{a-1}H_{n-m}\equiv H_{\lfloor p/3\rfloor}\pmod p$.\\
Now $p\equiv1\pmod3$, so by \cite[Lemma 17(2)]{lr}, we have $\binom{n}m\not\equiv0\pmod p$.
These, with (\ref{nm+1})-(\ref{pa-1n}) yield that we only need to prove that
\begin{align}\label{pa-1j}
p^{a-1}\sum_{\substack{j=0\\j\neq m}}^{n}\frac{\binom{n}j}{(j-m)(-4)^j}\equiv0\pmod p.
\end{align}
Now $n=(p^a-1)/2, m=(p^a-1)/6$. So by Fermat little theorem we have
$$
p^{a-1}\sum_{\substack{j=0\\j\neq m}}^{n}\frac{\binom{n}j}{(j-m)(-4)^j}\equiv-3(-1)^{(p^a-1)/2}p^{a-1}\sum_{\substack{j=0\\j\neq n-m}}^{n}\frac{\binom{n}j(-4)^j}{3j+1}\pmod{p}.
$$
There are only the items $3j+1=p^{a-1}(3k+1)$ with $k=0,1,\ldots,(p-1)/2$ and $k\neq (p-1)/3$, so by \cite[Theorem 1.2]{m-arxiv-2020} we have
\begin{align*}
&p^{a-1}\sum_{\substack{j=0\\j\neq m}}^{n}\frac{\binom{n}j}{(j-m)(-4)^j}\equiv-3(-1)^{\frac{p^a-1}2}\sum_{\substack{k=0\\k\neq (p-1)/3}}^{(p-1)/2}\frac{\binom{n}{kp^{a-1}+\frac{p^{a-1}-1}3}(-4)^{kp^{a-1}+\frac{p^{a-1}-1}3}}{3k+1}\\
&\equiv-3(-1)^{\frac{p^a-1}2}(-4)^{\frac{p^{a-1}-1}3}\binom{\frac{p^{a-1}-1}2}{\frac{p^{a-1}-1}3}\sum_{\substack{k=0\\k\neq (p-1)/3}}^{(p-1)/2}\frac{\binom{n}{k}(-4)^k}{3k+1}\equiv0\pmod{p}.
\end{align*}
Therefore the proof of Theorem \ref{Thsun} is complete.\qed

\vskip 3mm \noindent{\bf Acknowledgments.}
The author is funded by the Startup Foundation for Introducing Talent of Nanjing University of Information Science and Technology (2019r062).

\end{document}